\documentclass[a4paper,10pt]{article}
\usepackage[utf8]{inputenc}

\usepackage{a4wide}
\usepackage{amsmath}
\usepackage{amsfonts}
\usepackage{amssymb}
\usepackage{amsthm}
\usepackage{amstext}

\newtheorem{proposition}{Proposition}[section]
\newtheorem{corollary}[proposition]{Corollary}
\newtheorem{theorem}{Theorem}

\theoremstyle{definition}
\newtheorem{definition}[proposition]{Definition}
\newtheorem{remark}[proposition]{Remark}
\newtheorem{assumptions}[proposition]{Assumptions}

\usepackage{dsfont} 

\usepackage{xcolor}
\usepackage{tikz}
\usetikzlibrary{patterns}

\usepackage[numbers]{natbib}

\usepackage{mathtools}
\DeclarePairedDelimiter\floor{\lfloor}{\rfloor}

\definecolor{my_grey}{RGB}{170,170,170}   
\tikzset{
  schraffiert/.style={pattern=horizontal lines,pattern color=#1},
  schraffiert/.default=my_grey
}

\usepackage{microtype}

\usepackage{multirow}

\pdfoutput=1
\usepackage{hyperref}
\hypersetup{
    colorlinks=true,
    linkcolor=black,
    citecolor=black,
    filecolor=black,
    urlcolor=black,
}

\title{Method comparison with repeated measurements -- Passing-Bablok regression for grouped data with errors in both variables}
\author{\sc F. Baumdicker\footnote{Abteilung Mathematische Stochastik,
    Universit\"at Freiburg, D - 79104 Freiburg}
  $\mbox{}^,$\footnote{Corresponding author; email:
    baumdicker@stochastik.uni-freiburg.de} , U. Hölker$\mbox{}^*$
    }

\begin{document}

\date{}
\maketitle

\begin{abstract}


The Passing-Bablok and Theil-Sen regression are closely related non-parametric methods 
to estimate the regression coefficients and build tests on the relationship 
between the dependent and independent variables.
Both methods rely on the slopes of the connecting lines between pairwise measurements.
While Theil and Sen assume no measurement errors in the independent variable, 
the method from Passing and Bablok accounts for errors in both variables.
Here we consider the case where multiple, e.g. repeated, measurements with errors in both variables are available for $m$ samples.
We show that in this case excluding slopes between repeated measurements reduces the bias of the estimate.
We prove that the resulting Block-Passing-Bablok estimate for grouped data is asymptotically normally distributed.
If measurements of the independent variable are without error the variance of the estimate equals the variance of the Theil-Sen method with tied ranks.
If both variables are measured with imprecision the result depends on the fraction of measurements between groups
that fall within the range of each other. Only if no overlap between measurements of different groups occurs the variance
equals again the tied ranks version. Otherwise, the variance is smaller.
We explicitly compute this variance and provide a method comparison test 
for data with repeated measurements based on the method from Passing and Bablok for independent measurements.
If repeated measurements are considered this test has a higher power
to detect the true relationship between two methods.

\end{abstract}

keywords: Passing-Bablok regression, Theil-Sen regression, method comparison, method agreement, non-parametric statistics, rank-statistics,
          non-continuous sampling distribution, repeated observations, errors-in-variable, clinical chemistry

\section{Introduction}

Consider the simple linear regression, where $y= \alpha + \beta x$, for $x,y \in \mathbb R$.
Compared to ordinary least squares regression (OLS), non-parametric estimation of the regression coefficients $\alpha$ and $\beta$ can be more robust.
Deviations from normal distributed error terms, like heavy tails or outliers, 
do not disturb the estimation as strongly as for OLS. 
A good introduction to 
robust estimation is given in~\cite{wilcox2016introduction}.

Here we consider the non-parametric Theil-Sen regression (TSR) and Passing-Bablok regression (PBR),
which are both based on Kendall's rank correlation~\cite{Kendall1990rank} and provide a robust estimate of $\beta$.
PBR and TSR do not rely on the assumption of normally distributed errors.
The Theil-Sen estimate~\cite{Theil1950,Sen1968} for $\beta$ is given by the median of all slopes
of the connecting lines between pairwise measurements.
If measurement errors occur for both $x$ and $y$, the TSR is biased towards zero.
This phenomenon is also known from OLS and called regression dilution or attenuation.
For least squares, Deming regression instead of OLS can be used to account for errors in both variables.
In addition, \cite{Frost2000} suggest how to use repeated measurements to correct for regression dilution, if errors are normally distributed.
A variation of the Theil-Sen regression that accounts for errors in both variables is the Passing-Bablok estimator~\cite{Passing1983,Passing1984}.
The Passing-Bablok estimate is also given by the median of all slopes, but shifted by an offset $K$
to ensure that $x$ and $y$ are interchangeable. 
Interestingly, PBR and TSR seem to be popular in separate fields.
While TSR is popular in Metrology and Environmental Science,
PBR is mainly used in Clinical Biochemistry, Pharmacology and Laboratory Medicine to compare two alternative measurement methods.
This might be due to the fact that PBR is outlined in a guideline 
of the Clinical and Laboratory Standards Institute \cite{CLSIguideline}.
In \cite{Jensen2006} a protocol for method comparison studies in clinical laboratories is suggested, including a paragraph on PBR.
Furthermore, there have been attempts to use PBR for batch effect removal in gene expression analysis~\cite{Muller2016}.
Throughout this manuscript we will refer to the Passing-Bablok regression, 
but note that setting $K=0$ in the PBR equals the TSR and the presented results naturally hold for both methods.

We consider the case, where $p_k$ repeated measurements of two methods, for $m$ true pairs $(\tilde x_k,\tilde y_k)$ are available.
The alternative measurements of the two methods are given by data points $(x_{i,k},y_{i,k})$, $i \in \{1,\dots,p_k\}$.
This scenario includes simultaneous repeated measurements of the same sample with two methods,
as well as multiple measurements of the same sample separately measured for each method.
In the later, measurements need to be combined randomly into observation pairs $(x_{i,k}, y_{i,k})$.
If the data set for a PBR contains repeated measurements it is important to account for this.
In contrast to the expected slope $\beta$ between independent measurements, the expected slope between repeated measurements is $0$.
Hence, the slopes between points that correspond to repeated measurements of the same underlying true value
are meaningless and would distort PB estimates and lower the power of the associated statistical test if included.
Here we describe how the variance changes if we omit the meaningless slopes between repeated measurements and provide the resulting
test on the equivalence of two methods.

The paper is organized as follows:
First, we consider the assumptions and definitions 
of Passing and Bablok and introduce the necessary adaptation for repeated measurements in section~\ref{section2}.
Section~\ref{section3} states our results, 
including the asymptotic confidence interval for the estimated parameter $\hat \beta$ in Corollary~\ref{T3}.
The test for the equivalence of two methods with repeated measurements is given in Remark~\ref{cor:test}.
In section~\ref{section4} we discuss the implications of the suggested Block-Passing-Bablok procedure and compare it empirically to the PBR without repeated measurements.
Finally, proofs are given in section~\ref{section5}.

\section{The Passing-Bablok regression for repeated measurements}

\label{section2}

\subsection{The standard Passing-Bablok regression for independent measurements}

Under the hypotheses of a structural linear relationship between two measurement methods, i.e.
\begin{align*}
	\tilde y_i = \alpha + \beta \tilde x_i,
\end{align*}
for $ \tilde x_i, \tilde y_i \in \mathbb R$, where $y$ corresponds to one and $x$ to the other method, 
Passing and Bablok~\cite{Passing1983} considered the following assumptions:

\begin{assumptions}{(standard Passing-Bablok Regression)}
	\begin{itemize}
	\item[i)]All points $(x_{i},y_{i})$ are of the form
		\begin{align*}
			x_{i} = \tilde{x}_i + \varepsilon_{i}, \\
			y_{i} = \tilde{y}_i + \eta_{i},
		\end{align*}
		with $\tilde x_{i},\tilde y_{i}$ the (non-random) 'true' values of the measurements and
		error terms $\varepsilon_{i}, \eta_{i}$ and $i \in \{1, \dots, n\}$.
	\item[ii)] $\varepsilon_1,\dots,\varepsilon_n$ are iid and come from an arbitrary continuous distribution with mean zero, and \label{VssPB}
	          $\eta_1,\dots,\eta_n$ are iid and come {\color{black} from the same type of distribution such that $\beta \varepsilon_i$ and $\eta_i$ are iid.}
	          Moreover, $\varepsilon_i$ and $\eta_j$ are independent for each $i,j$.
	\end{itemize}
\end{assumptions}

\subsection{Block-Passing-Bablok regression for multiple dependent measurements}

Sometimes the data at hand is not independent as stated in assumption~\ref{VssPB}(ii) but grouped into dependent sets of measurements.
Examples include common situations, e.g. if measurements have been repeated multiple times
for the same sample or a series of measurements is done under the same conditions.
In this case, the Passing-Bablok Regression cannot be used straight away.
Here, we expand the Passing-Bablok Regression such that it considers the data to be available in $m$ groups with $p_k (k = 1,\dots,m)$ members.
We index our data points $(x_{i,k},y_{i,k})$ by group $k \in \{1, \dots m\}$ and individual $i \in \{1, \dots p_k\}$ from this group.
Each $(x_{i,k},y_{i,k})$ represents a measurement of the true values $(\tilde x_k, \tilde y_k)$ and
we assume again a linear relationship
\begin{align*}
	\tilde y_{k} = \alpha + \beta \tilde x_{k}.
\end{align*}
We make the following assumptions for repeated measurement data:

\begin{assumptions}{(Passing-Bablok Regression for grouped data)}
\label{VssPBgrouped}
	\begin{itemize}
	\item[i)]All $n=\sum_{k=1}^{m} p_k $ points $(x_{i,k},y_{i,k})$ are of the form
		\begin{align*}
			x_{i,k} = \tilde{x}_k + \varepsilon_{i,k}, \\
			y_{i,k} = \tilde{y}_k + \eta_{i,k},
		\end{align*}
		with $\tilde{x}_k, \tilde{y}_k$ the 'true' values of the measurements in group $k$ and 
		error terms $\varepsilon_{i,k}, \eta_{i,k}$ indexed by group $k \in \{1, \dots m\}$ and individual $i \in \{1, \dots, p_k\}$ from this group. 
	\item[ii)] All $\varepsilon_{i,k}$ are iid and come from an arbitrary and continuous distribution with mean zero, and
	           all $\eta_{i,k}$ are iid and come  {\color{black} from the same type of distribution such that $\beta \varepsilon_{i,k}$ and $\eta_{i,k}$ are iid.}
	           Moreover, $\varepsilon_{i,k}$ and $\eta_{j,l}$ are independent for each $i,j,k,l$.
	\end{itemize}
\end{assumptions}

\begin{assumptions}{(non overlapping groups)}
	  \label{VssPBnooverlap}
         \\All $m$ groups are strictly separated on the x-axis, i.e. almost surely\ 
	 \begin{align*}
	  \max_{1 \leq i \leq p_k} x_{i,k} < \min_{1 \leq j \leq p_l} x_{j,l} \quad \forall ~ k,l \in \{1, \dots m\}, k < l,
	 \end{align*}			
\end{assumptions}

\begin{remark}
 Assumption~\ref{VssPBnooverlap} can only be fulfilled if the error distribution is not supported on the whole real line, e.g.\ for uniformly distributed errors
 if $| \tilde x_i - \tilde x_j|$ is large enough for all $i \neq j$.
 Note further that groups need to be also separated at the $y$-axis as in Assumption~\ref{VssPBnooverlap} to maintain the symmetry of the method.
 However, we will mainly consider the general case with possibly overlapping groups such that only Assumptions~\ref{VssPBgrouped}~i)-ii) hold.
\end{remark}

{\color{black}
\begin{remark}
 Assumption~\ref{VssPB} ii) and~\ref{VssPBgrouped} ii) set the measurement error ratio equal to the regression slope parameter $\beta$, i.e. for the variances of error terms, 
 $\frac{\sigma_\eta^2}{\sigma_\varepsilon^2} = \beta^2$ holds.
 The assumption ensures that the slopes are symmetrically distributed around $\beta$ (see remark~\ref{slopesymmetric}).
 A common case where the assumptions are fulfilled is the comparison of two measurement methods
 under the null hypothesis that both methods are equivalent, i.e. $\beta = 1$, and all $\epsilon_i$ and $\eta_i$ are iid.
\end{remark}
}

\subsection*{Estimating the regression parameter $\beta$ from grouped data}

The Passing-Bablok regression for independent measurements~\cite{Passing1983}
makes use of the slopes $S_{i,j} = \frac{y_{i} - y_{j}}{x_{i} - x_{j}} $
between all pairs of points $(x_i,y_i)$ and $(x_j,y_j)$.
Under Assumptions~\ref{VssPB}, it is easy to see that

\begin{equation}
\mathbb E[S_{i,j}] = \beta (1 - \mathbb E[\frac{\epsilon_i - \epsilon_j}{ \tilde x_i - \tilde x_j +  \epsilon_i - \epsilon_j} ] )  \label{expectedslope}
\end{equation}
which is closer to $\beta$ the larger $|\tilde x_i - \tilde x_j|$.
In contrast, under Assumptions~\ref{VssPBgrouped} we get that the expected slope between $y_{i,k}$ and $y_{j,l}$ 
equals \eqref{expectedslope} if $k\neq l$ and is zero for $k=l$.
A similar argument holds for the median of slopes between and within groups.
To include only the meaningful slopes for the estimation of the regression parameter~$\beta$,
we have to exclude all repeated measurement pairs within each group from the analysis.

\begin{definition}[Block-Passing-Bablok regression for grouped data]%
\label{def:extendedPB}
\quad
\begin{enumerate}
 \item For the estimation of the regression parameters $\alpha$ and $\beta$, 
 compute the slopes of the connecting lines between any pair of points from different groups.
 The slopes are given by
 \begin{align*}
	 S_{ij}^{kl} &= \frac{y_{i,k} - y_{j,l}}{x_{i,k} - x_{j,l}} 
	 \qquad \text{for } k \neq l \text{, } i= 1, \dots , p_k \text{ and } j= 1, \dots , p_l .
 \end{align*}
\item Discard identical measurements as well as all slopes with $S_{ij}^{kl} = -1$.
Thereby obtain $N \leq \binom{n}{2} - \sum\limits_{k=1}^m \binom{p_k}{2}$ slopes.

\item The slope parameter $\beta$ is estimated by the shifted median of the slopes with the offset
      \begin{align}
	K = \# \{S_{ij}^{kl} | S_{ij}^{kl}<-1\}.
	\label{valueK}
      \end{align}
      I.e. if the ranked sequence of slopes is given by $S_{(1)} \leq S_{(2)} \leq \dots \leq S_{(N)}$,
      $\beta$ is estimated by 
\begin{align*}
 \hat{\beta} =
		\begin{cases}
			S_{(\frac{N+1}{2}+K)} ,& \text{if $n$ is odd}, \\
			\frac{1}{2} \big( S_{(\frac{N}{2} + K)} + S_{(\frac{N}{2} + K +1)}    \big) ,& \text{if $n$ is even}
		\end{cases}	
\end{align*}
\item The intercept parameter $\alpha$ is estimated by
\begin{equation*}
 \hat{\alpha} = \text{median} \{ y_i -\hat{\beta}x_i \}
\end{equation*}
\end{enumerate}
\end{definition}

\begin{remark}
Naturally, the description of the groupwise method in Definition~\ref{def:extendedPB} 
is close to the description of the original method introduced by Passing and Bablok.
If we set $m = n$ and $p_k=1$ we regain the classical estimator for the regression parameter $\beta$~\cite{Passing1983}.
In this case, we compute the slopes of the connecting lines between any pair of points, that are given by
\begin{align*}
	S_{11}^{kl}  = \frac{y_{1,k} - y_{1,l}}{x_{1,k} - x_{1,l}} \qquad \text{for} \quad 1 \leq k < l \leq n.
\end{align*}
Without considering groups there are $\binom{n}{2}$ possible lines to connect any two points of an $n$-dimensional data set.
Again two identical measurements with $ x_i = x_j$ and $y_i = y_j$ are not considered for the estimation.
Further, any slopes with a value of $-1$ are disregarded, such that
we have at most $\binom{n}{2}$ slopes to consider.
\end{remark}

\begin{remark}
Note that setting $K=0$ in Definition~\ref{def:extendedPB} results in the Theil-Sen estimator.
To get an estimator where methods can be used interchangeably, Passing and Bablok defined the offset determined by $K$.
The definition of $K$ as the number of slopes with a value smaller than $-1$ in equation~\eqref{valueK} corresponds 
to the null hypothesis $\beta = 1$, and needs to be adapted for other hypothesis for the value of $\beta$.
If the null hypothesis is not true, setting $K$ in this way introduces a bias towards higher estimates for $\beta$.
The median slope between independent measurements within one group with offset $K$ is no longer given by zero if $K>0$.
This effect can be seen in Table~\ref{tab:performance}, e.g. for $\beta=0.2$ for both the classic and the Block-PBR.
Since the offset $K$ drives the median of meaningless slopes between repeated measurements towards $1$, this can also lead 
to overconfidence for $\beta = 1$ if the classic PBR is used. Table~\ref{tab:performance} for $\beta=0.8$ illustrates this effect.
\end{remark}

\section{Results}

\label{section3}

\begin{definition}
 \label{defCtilde}
 Under Assumptions~\ref{VssPBgrouped} let $\tilde C$ be the difference between the number of slopes 
 that are larger than $\beta$ and the number of slopes lower than $\beta$.
 I.e. for 
 \begin{align*}
	P(\beta) =& ~ \# \{ ((i,k)(j,l)) ~|~ S_{ij}^{kl} > \beta \}\text{ and} \\
	Q(\beta) =& ~ \# \{ ((i,k)(j,l)) ~|~ S_{ij}^{kl} < \beta \}\text{ we define }
\end{align*}
\begin{equation}
\tilde C :=  P(\beta) - Q(\beta).
\end{equation}
 Note that we will denote this number by $C$ instead of $\tilde C$ in the special case of Assumptions~\ref{VssPB} where no groups are considered.
\end{definition}

{\color{black}
\begin{remark}
\label{slopesymmetric}
Using assumption~\ref{VssPBgrouped} ii) we can write 
\begin{equation*}
 S_{i,j}^{k,l} = \beta \frac{x_{i,k} - x_{j,l} + z_{i,j}^{k,l} }{x_{i,k} - x_{j,l} + z_{i,j}^{'k,l}}
\end{equation*}
where $z_{i,j}^{k,l} = (\eta_{i,k} - \eta_{j,l}) \beta^{-1}$ and $z_{i,j}^{'k,l} = \epsilon_{i,k} - \epsilon_{j,l}$ are independent and from the same distribution~\cite{Passing1983}.
Thus the number of slopes above and below $\beta$ and consequently $\tilde C$ are free from the regression parameters.
\end{remark}
}

\begin{theorem}[Variance of $\tilde C$]
\label{T1}
If Assumptions~\ref{VssPBgrouped} (i)-(ii) hold and 
$q_{ku}$ is the expected fraction of triplets with one point from group $u$ and two points from group $k$ with 
$\min(x_{i,k},x_{j,k}) < x_{s,u} < \max(x_{i,k},x_{j,k})$ and\\

(a) if the group sizes are given by $(p_1,\dots,p_m)$,
\begin{equation}
  \mathbb V[\tilde C] = \frac{1}{18} \bigg( n(n-1)(2n+5)  - \sum\limits_{k=1}^m p_k (p_k-1)((2p_k+5) + 4 \sum\limits_{u\neq k}^{m} p_u q_{ku} )   \bigg). \label{VarC}
\end{equation}

(b) Consequently, if the group sizes are equal, i.e.\ $p_k = p := \tfrac{n}{m}$ for $k = 1,\dots,m$,
\begin{equation}
\mathbb V[\tilde C] = \frac{n}{18} \big( 3(n-p) + 2 (n^2 - p^2)  \big) 
   - \frac29 p^2 (p-1) \sum\limits_{k=1}^m \sum\limits_{u\neq k}^m q_{ku}.
   \label{VarCequal} 
\end{equation}

\end{theorem}

\begin{remark}[classical Passing-Bablok regression]\quad\\
If we set $p=1$ and $m=n$ in Theorem~\ref{T1} we regain the result from Passing and Bablok~\cite{Passing1983} and Theil-Sen~\cite{Sen1968}
for independent measurements, where
$$ \mathbb V[C] = \frac{n(n-1)(2n+5)}{18}.$$
\end{remark}

\begin{remark}
In Theorem~\ref{T1}, note that (b) is a direct consequence from (a) as setting $p_k = p$ in equation~\eqref{VarC} gives
\begin{equation*}
 \frac{1}{18} \bigg( n(n-1)(2n+5)  - mp (p-1)((2p+5) - 4 p^2(1-p)\sum\limits_{k=1}^{m}\sum\limits_{u\neq k}^{m} q_{ku} )   \bigg),
\end{equation*}
which leads to \eqref{VarCequal} as $mp=n$.
\end{remark}

\begin{theorem}[asymptotic normality of $\tilde C$]
\label{T2}\quad\\
Let the number of groups $m$ be fixed and consider the limit of a large samplesize $n = \sum\limits_{i=1}^m p_k$.

(a) If the group sizes are given by $(p_1(n),\dots,p_{m}(n))$ and Assumptions~\ref{VssPBgrouped} (i)-(ii) hold,
$\tilde C$ is asymptotically normally distributed with mean zero and variance given by
$$\mathbb V[\tilde C] \sim \frac19 n^3 \Big(
1 - l_m - l_o 
\Big)\text{, with}$$
\begin{align*}
 l_m &= \lim_{n\to\infty}\sum\limits_{k=1}^m \frac {p_k(n)^3}{n^3}\text{,}\\
 l_o &= \lim\limits_{n\to \infty} \sum\limits_{k=1}^m \sum\limits_{u\neq k}^m n^{-3} \big(p_k(n)\big)^2 p_u(n) q_{ku}\text{,} 
\end{align*}
where $q_{ku}$ is the expected fraction of triplets with one point from group $u$ and two points from group $k$ with 
$\min(x_{i,k},x_{j,k}) < x_{s,u} < \max(x_{i,k},x_{j,k})$.


(b) More specifically, for $m$ groups with equal group sizes as in Theorem~\ref{T1}(b) and if Assumptions~\ref{VssPBgrouped} and~\ref{VssPBnooverlap} hold,
$\tilde C$ is asymptotically normally distributed with mean zero and variance 
$$\mathbb V[\tilde C] \sim \frac19 n^3 (1 - \frac {1}{m^2} ). $$

\end{theorem}

\begin{corollary}[confidence interval for $\beta$]
\label{T3}\quad\\
Let $w_{\frac{\gamma}{2}}$ denote the $\big( 1- \frac{\gamma}{2} \big)$-quantile of the standardized normal distribution.
Let further $\tilde \sigma := (\mathbb V[\tilde C])^{\frac 12}$, 
\begin{equation*}
	\tilde C_\gamma = w_{\frac{\gamma}{2}}\tilde\sigma,\quad 
	M_1 = ~ \floor[\Big]{\frac{N - \tilde C_\gamma}{2}}, \text{ and}\quad
	M_2 = ~N -M_1 + 1.
\end{equation*}
The asymptotic confidence interval for $\beta$ with significance level $1-\gamma$ is then given by 
\begin{align}
	I = \big[ S_{(M_1 + K)}, S_{(M_2 + K)} \big].
	\label{KIbeta}
\end{align}	  
\end{corollary}

\begin{remark}[confidence interval for $\alpha$]\quad\\
 Let $b_L$ denote the lower and $b_U$ denote the upper limit of the confidence interval $I$ for the slope $\beta$ in equation~\eqref{KIbeta}. Then,
\begin{align}
	a_L =& ~ \text{median} \{y_i - b_U x_i \}, \quad \text{and} \qquad \\
	a_U =& ~ \text{median} \{y_i - b_L x_i \}
	\label{KIalpha}
\end{align}
are the corresponding limits for the intercept $\alpha$. 
\end{remark}

\begin{remark}[Equivalence of two methods with repeated measurements]
\label{cor:test}
Given two measurement methods with measurements given by $x_{i,k}$ and $y_{i,k}$.
\begin{itemize}
 \item The equivalence of both methods can be concluded if $0\in[a_L,a_U]$ and $1\in[b_L, b_U]$. 
 \item An intercept $\hat \alpha \neq 0$ indicates a constant systemic difference between the two measurement methods.
 \item A regression parameter $\hat \beta \neq 1$ indicates a proportional systemic difference between the measurement methods.
\end{itemize}
Note that here the asymptotic confidence level is calculated for $\beta$ and not $\alpha$.
\end{remark}

\begin{remark}[non-overlapping groups]
In case of a data set where Assumption~\ref{VssPBnooverlap} is fulfilled, i.e.\ with non-overlapping groups on the x-axis, 
we can set $q_{ku} = 0$ and formula (\ref{VarC}) transforms to
\begin{equation}
  \mathbb V[\tilde C] = \frac{1}{18} \bigg( n(n-1)(2n+5)  - \sum\limits_{k=1}^m p_k (p_k-1)(2p_k+5)   \bigg). \label{VarCnonoverlapping}
\end{equation}
And consequently for equal group sizes
\begin{equation}
\mathbb V[\tilde C] = \frac{n}{18} \big( 3(n-p) + 2 (n^2 - p^2)  \big). \label{VarCequalnonoverlapping} 
\end{equation}

Note that \eqref{VarCnonoverlapping} equals the correction for tied ranks in one variable~\cite{Sen1968}.
To see this consider the fact that setting $x_{i,k}$ to the mean of the corresponding group for $k=1,\dots,m$ 
does not change the sign of slopes between separated groups.
So compared to the variance for non-overlapping groups the variance for overlapping groups is always smaller.
For overlapping groups the test for the hypotheses $\beta = 1$ from Remark~\ref{cor:test} based on equation~\eqref{VarCnonoverlapping}
is thus a conservative test for the equivalence of the measurement methods,
but the power of the test could in principle be improved if a reliable estimate for $q_{ku}$ is available.
\end{remark}

\section{Discussion and empirical comparisons of \\regular and Block-Passing-Bablok Regression}

\label{section4}

Passing-Bablok Regression is recommended as a robust method such that 
extreme values can be included and the errors do not have to be normally distributed~\cite{Jensen2006}.
However, the assumption of independent measurements in the classical Passing-Bablok regression seems to be frequently not fulfilled.
Duplicated and repeated measurements are often an important part of studies to assess the variance of measurements and help to identify outliers \cite{Jensen2006}.
As a random example among various studies including repeated measurements consider Figure 2 in~\cite{Weber695}, where measurements of different patients have been repeated in different numbers.
But to which degree does this harm the results of a classical Passing-Bablok regression?

We recall, the Block-Passing-Bablok Regression only considers slopes
between pairs of points from different groups for estimating the regression parameter $\beta$.
Since the regression parameters are estimated as medians, a relatively small difference for the number of considered slopes
does not have a large influence on the estimations.
Hence if the group sizes are equal, instead of $\binom{n}{2}$ slopes as in the original method, we utilize approximately $\binom{n}{2}(1-\frac 1m)$ slopes.
As shown in Theorem~\ref{T2}~(a) the discrepancy between the unadapted and the Block-Passing-Bablok Regression
will quickly vanish for a large number of separated and equally sized groups.

In a hypothetical example with varying group sizes where $p_1 \gg p_2,\dots,p_m$, the Block-Passing-Bablok Regression yields better results than the original method.
The median of all slopes between a very large and a very small group would be a slope within the large group.
Since we assume random errors within groups, this median slope would be sampled from a set of meaningless slopes with mean $0$,
which transforms into $1$ if the offset $K$ is used. This biases the estimate of $\beta$ towards $1$ and lowers the power to detect $\beta \neq 1$.
The Block-Passing-Bablok Regression provides more reliable estimates in this setting.

\begin{figure}

\includegraphics[width=0.49\textwidth]{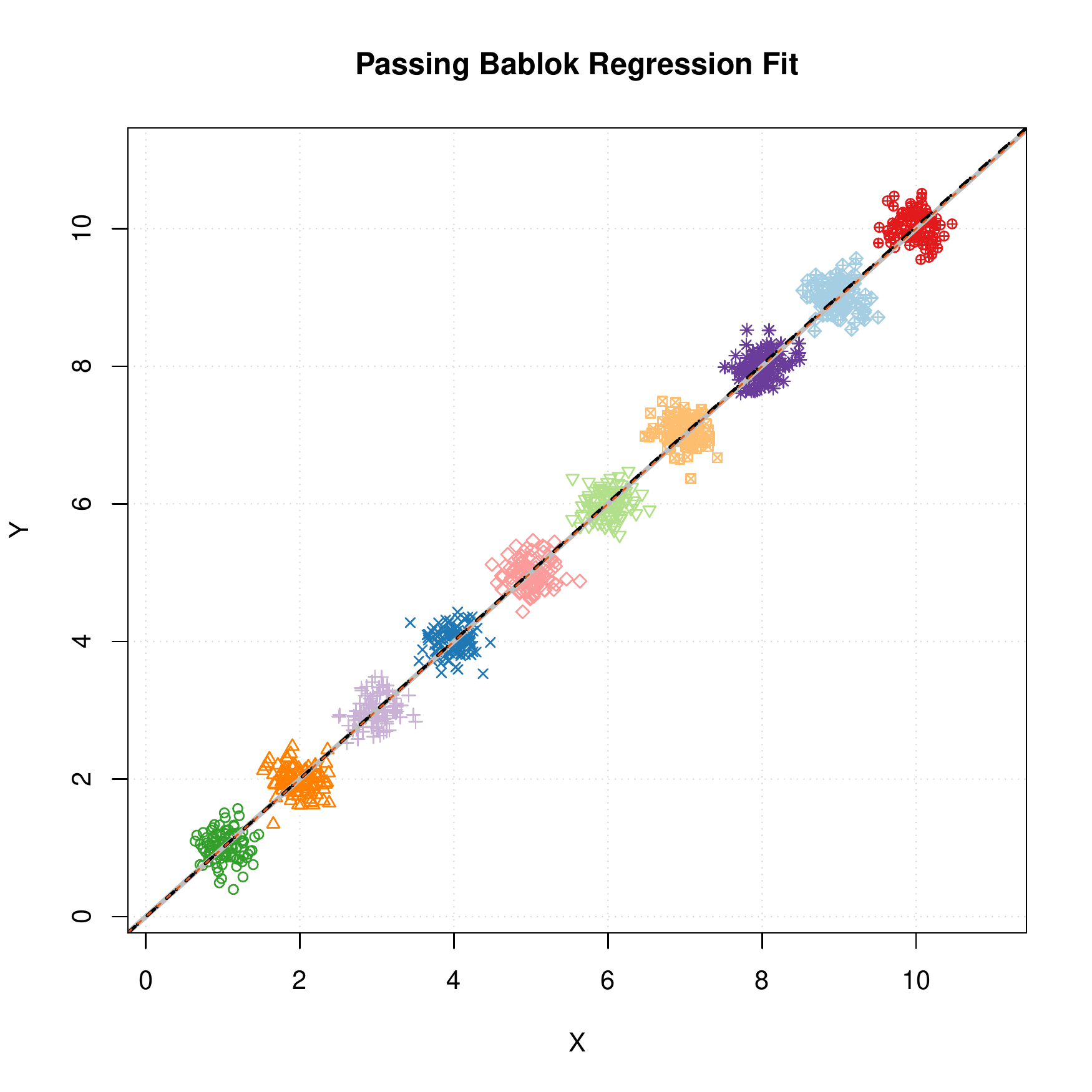}
\includegraphics[width=0.49\textwidth]{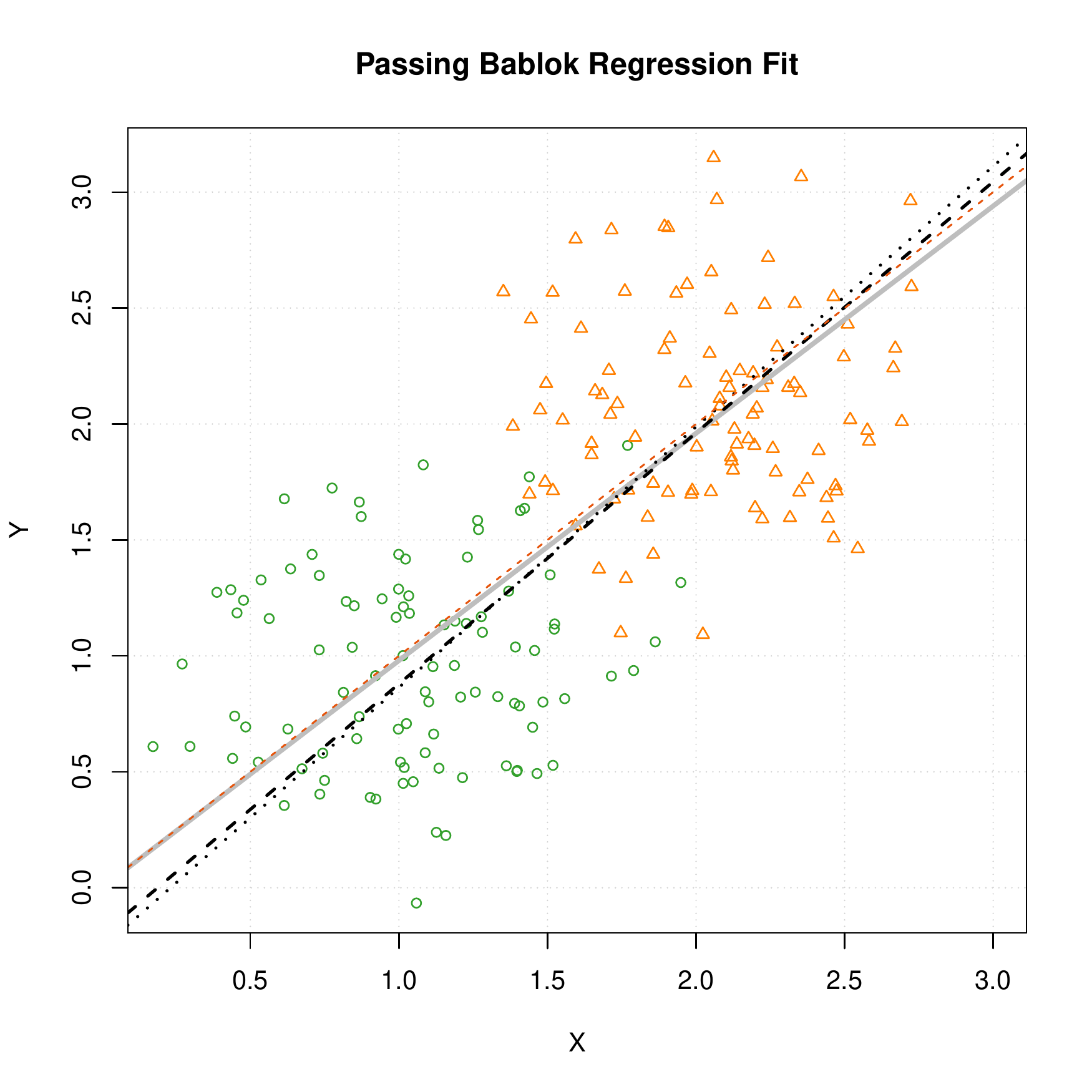}

\includegraphics[width=0.49\textwidth]{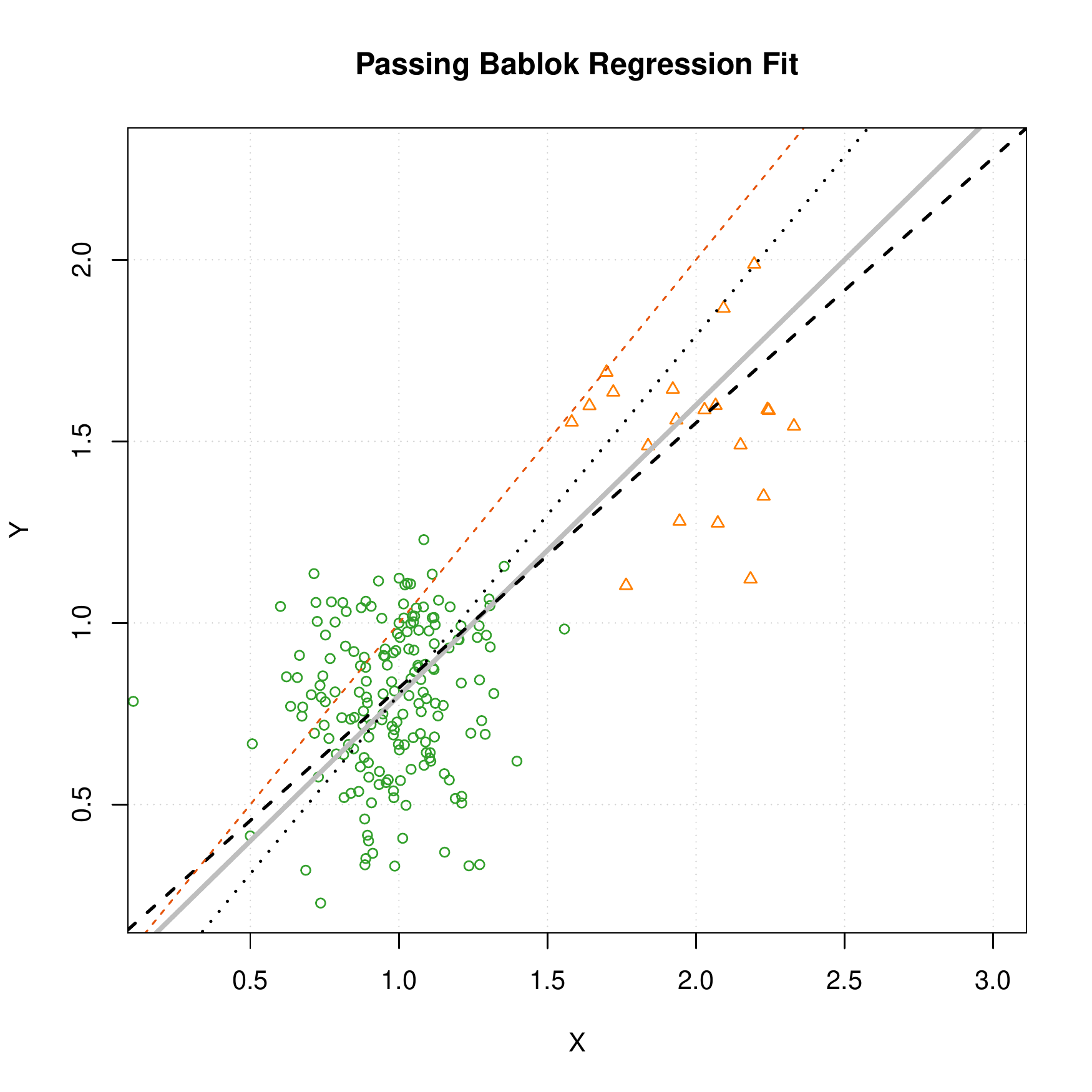}
\includegraphics[width=0.49\textwidth]{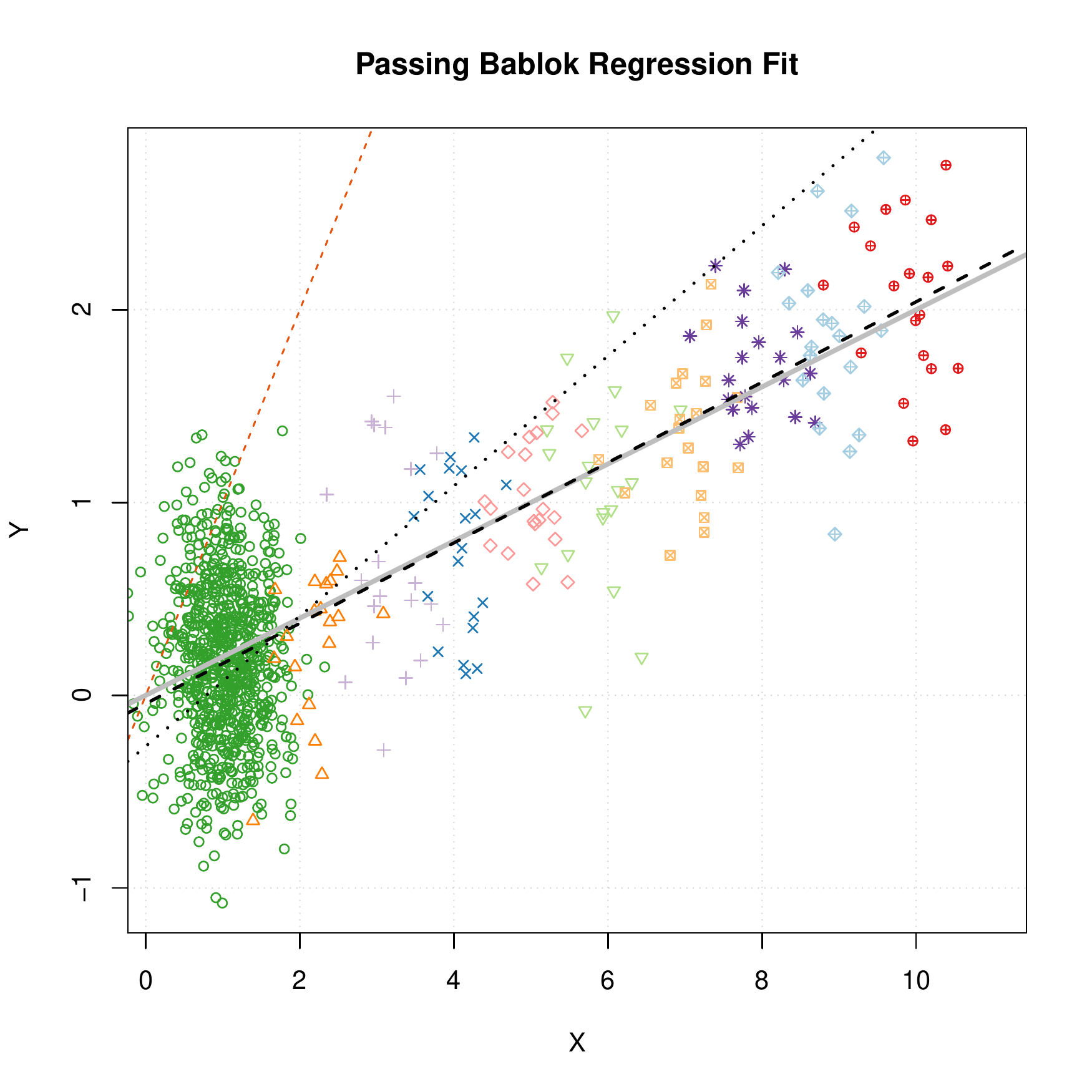}

\caption{
Four scenarios of our simulations are shown. The solid gray line represents the true slope $\beta$. The black long dashed line is the Block-Passing-Bablok estimate
and the dotted line is the classic Passing-Bablok estimate ignoring the group memberships.
Group memberships are illustrated via the shape and color of the measurements $(x_{i,k},y_{i,k})$. The red short dashed line shows the identity.
Top left: $\beta = 1, p_1=\dots=p_{10}=100, \sigma_\epsilon = \sigma_\eta = 0.2 $;
Top right: $\beta = 0.98, p_1=p_2=100, \sigma_\epsilon = \sigma_\eta = 0.4 $;
Bottom left: $\beta = 0.8, p_1=180,p_2=20, \sigma_\epsilon = \sigma_\eta = 0.2 $;
Bottom right: $\beta = 0.2, p_1=820,p_2=\dots=p_{10}=20, \sigma_\epsilon = \sigma_\eta = 0.4 $.
}
\label{paragridsample}
\end{figure}

To illustrate the influence of the group sizes and the overlap between groups on the estimation of $\beta$ and the power of the test of the hypotheses $\beta = 1$ in Remark~\ref{cor:test}
we simulated 16 illustrative scenarios. Evaluations of the Passing-Bablok Estimates have been performed with the \textit{mcr} package \cite{mcrpackage}
for the classic PBR and with a custom script based on \textit{mcr} for the Block-PBR. The script is available from the authors upon request.
A subset of the considered settings is illustrated in Figure~\ref{paragridsample}.
The results of the simulations are shown in Table~\ref{tab:performance}. 

\begin{table}[h]
\centering 
{\small
        \begin{tabular}{lllrrccrrrr}
            \multirow{2}{*}{slope} & \multirow{2}{*}{groups} & \multirow{2}{*}{overlap} & \multicolumn{2}{c}{mean $\hat b$} & \multicolumn{2}{c}{mean $I$} & \multicolumn{2}{c}{$\mathbb P(\beta \in I)$} & \multicolumn{2}{c}{$\mathbb P(1 \notin I)$}  \\
             &  &  & cPBR & BPBR & cPBR & BPBR & cPBR & BPBR & cPBR & BPBR \\ 
            \hline\hline
            \multirow{8}{*}{$\beta=1.0$} & \multirow{2}{*}{100-100} & \multicolumn{1}{l}{low}
				& 1.001 & 1.001 & [.923,1.085] & [.923,1.085] & .950 & .950 & .050 & .050    
				 \\
                                 &  & \multicolumn{1}{l}{high}
                                 & 1.003 & 1.004 & [.870,1.156] & [.856,1.177] & .950 & .949 & .050 & .051
                                 \\\cline{2-11}
                                 & \multirow{2}{*}{180-20} & \multicolumn{1}{l}{low} 
                                 & 1.003 & 1.002 & [.860,1.169] & [.874,1.148] & .951 & .949 & .049 & .051
                                 \\
                                 &  & \multicolumn{1}{l}{high} 
                                 & 1.005 & 1.01 & [.812,1.246] & [.771,1.326] & .951 & .947 & .049 & .053
                                 \\\cline{2-11}
                                 & \multirow{2}{*}{10$\times$100} & \multicolumn{1}{l}{low} 
                                 & 1 & 1 & [.994,1.006] & [.994,1.006] & .950 & .950 & .050 & .050
                                 \\
                                 &  & \multicolumn{1}{l}{high} 
                                 & 1 & 1 & [.988,1.012] & [.987,1.012] & .952 & .952 & .048 & .048
                                 \\\cline{2-11}
                                 & \multirow{2}{*}{820-9$\times$20} & \multicolumn{1}{l}{low} 
                                 & 1 & 1 & [.988,1.013] & [.991,1.009] & .951 & .952 & .049 & .048
                                 \\
                                 &  & \multicolumn{1}{l}{high} 
                                 & 1 & 1 & [.977,1.024] & [.982,1.018] & .951 & .951 & .049 & .049
                                 \\\hline
            \multirow{8}{*}{$\beta=0.98$} & \multirow{2}{*}{100-100} & \multicolumn{1}{l}{low}
				 & .983 & .981 & [.907,1.067] & [.904,1.064] & .950 & .950 & .071 & .076
				 \\
                                 &  & \multicolumn{1}{l}{high} 
                                 & .989 & .984 & [.856,1.141] & [.838,1.156] & .950 & .948 & .054 & .057
                                 \\\cline{2-11}
                                 & \multirow{2}{*}{180-20} & \multicolumn{1}{l}{low} 
                                 & .991 & .983 & [.850,1.157] & [.856,1.128] & .949 & .947 & .053 & .061
                                 \\
                                 &  & \multicolumn{1}{l}{high} 
                                 & .998 & .991 & [.804,1.239] & [.753,1.306] & .948 & .951 & .049 & .051
                                 \\\cline{2-11}
                                 & \multirow{2}{*}{10$\times$100} & \multicolumn{1}{l}{low} 
                                 & .980 & .980 & [.974,.986] & [.974,.986] & .950 & .949 & 1 & 1
                                 \\
                                 &  & \multicolumn{1}{l}{high} 
                                 & .980 & .980 & [.968,.993] & [.968,.993] & .950 & .951 & .876 & .880
                                 \\\cline{2-11}
                                 & \multirow{2}{*}{820-9$\times$20} & \multicolumn{1}{l}{low} 
                                 & .981 & .980 & [.969,.994] & [.971,.989] & .950 & .951 & .838 & .994
                                 \\
                                 &  & \multicolumn{1}{l}{high} 
                                 & .982 & .980 & [.959,1.006] & [.963,.998] & .946 & .948 & .326 & .607
                                 \\\hline
            \multirow{8}{*}{$\beta=0.8$} & \multirow{2}{*}{100-100} & \multicolumn{1}{l}{low} 
				 & .828 & .801 & [.757,.906] & [.730,.877] & .888 & .953 & .987 & .998
				 \\
                                 &  & \multicolumn{1}{l}{high} 
                                 & .854 & .807 & [.730,.997] & [.672,.963] & .881 & .952 & .536 & .679
                                 \\\cline{2-11}
                                 & \multirow{2}{*}{180-20} & \multicolumn{1}{l}{low} 
                                 & .888 & .802 & [.754,1.051] & [.686,.934] & .769 & .949 & .303 & .823
                                 \\
                                 &  & \multicolumn{1}{l}{high} 
                                 & .924 & .812 & [.738,1.158] & [.592,1.095] & .770 & .950 & .117 & .305
                                 \\\cline{2-11}
                                 & \multirow{2}{*}{10$\times$100} & \multicolumn{1}{l}{low} 
                                 & .801 & .800 & [.795,.807] & [.794,.806] & .934 & .948 & 1 & 1
                                 \\
                                 &  & \multicolumn{1}{l}{high} 
                                 & .802 & .800 & [.791,.813] & [.789,.812] & .935 & .947 & 1 & 1
                                 \\\cline{2-11}
                                 & \multirow{2}{*}{820-9$\times$20} & \multicolumn{1}{l}{low} 
                                 & .813 & .800 & [.802,.825] & [.792,.808] & .398 & .953 & 1 & 1
                                 \\
                                 &  & \multicolumn{1}{l}{high} 
                                 & .824 & .800 & [.803,.847] & [.784,.816] & .385 & .948 & 1 & 1
                                 \\\hline
            \multirow{8}{*}{$\beta=0.2$} & \multirow{2}{*}{100-100} & \multicolumn{1}{l}{low} 
				 & .317 & .202 & [.257,.383] & [.144,.261] & .022 & .948 & 1 & 1
				 \\
                                 &  & \multicolumn{1}{l}{high} 
                                 & .461 & .279 & [.352,.588] & [.165,.404] & .001 & .758 & 1 & 1
                                 \\\cline{2-11}
                                 & \multirow{2}{*}{180-20} & \multicolumn{1}{l}{low} 
                                 & .559 & .201 & [.415,.753] & [.105,.302] & 0 & .952 & .975 & 1
                                 \\
                                 &  & \multicolumn{1}{l}{high} 
                                 & .508 & .280 & [.318,.733] & [.093,.501] & 0 & .906 & .685 & 1
                                 \\\cline{2-11}
                                 & \multirow{2}{*}{10$\times$100} & \multicolumn{1}{l}{low} 
                                 & .204 & .200 & [.200,.209] & [.196,.205] & .567 & .954 & 1 & 1
                                 \\
                                 &  & \multicolumn{1}{l}{high} 
                                 & .212 & .204 & [.203,.221] & [.195,.213] & .236 & .846 & 1 & 1
                                 \\\cline{2-11}
                                 & \multirow{2}{*}{820-9$\times$20} & \multicolumn{1}{l}{low} 
                                 & .274 & .200 & [.255,.302] & [.194,.206] & 0 & .951 & 1 & 1
                                 \\
                                 &  & \multicolumn{1}{l}{high} 
                                 & .328 & .202 & [.299,.366] & [.189,.215] & 0 & .941 & 1 & 1
                                 \\\hline
        \end{tabular}
}
\caption{\label{tab:performance} We consider four different combinations of group sizes denoted by 100-100 $(p_1 = p_2  =100)$, 180-20 $(p_1 = 180, p_2  =20)$,
10$\times$100 $(p_1 = \dots = p_{10}  =100)$, and 820-9$\times$20 $(p_1 = 820, p_2 = \dots = p_{10}  =20)$.
The true values of the samples are given by $\tilde x_k = k$ and $\tilde y_k = b \tilde x_k$, for $k = 1, \dots,m$.
In addition, we vary the variance within groups, setting $\sigma_\epsilon = \sigma_\eta$ to either $0.2$ or $0.4$ to create high or low overlap between groups.
We simulate the performance under the null hypothesis $\beta=1$ and alternative scenarios where $\beta \in \{0.98,0.8,0.2\}$.
The difference between the classical Passing-Bablok regression (cPBR) and the 
Block-Passing-Bablok regression (BPBR) from Definiton~\ref{def:extendedPB} is shown.
In particular, we provide estimations for $\mathbb E[\hat \beta]$, $\mathbb E[I]$, 
the probability that the true slope $\beta \in I$ and the probability to reject the null hypothesis $\beta = 1$.
} 
\end{table}

\section{Proofs}

\label{section5}

For Theorem~\ref{T1}, Theorem~\ref{T2}, and the verification of confidence intervals in Corollary~\ref{T3} we need to compute the variance of 
\begin{align*}
	\tilde{C} : = \# \{ S_{ij}^{kl} ~|~ S_{ij}^{kl} > \beta \} -  \# \{ S_{ij}^{kl} ~|~ S_{ij}^{kl} < \beta \}
\end{align*}
the number of slopes $S_{ij}^{kl}$ bigger than $\beta$ minus the number of slopes $S_{ij}^{kl}$ smaller than $\beta$, 
as introduced in Definition~\ref{defCtilde}.\\

To simplify the notation for the proofs, we will instead consider a transformed dataset, where 
$x_{i,k}$ are divided by $\beta$ such that the errors in both dimensions are identically distributed (see Assumption~\ref{VssPBgrouped}(ii)).
If we also substract $\beta^{-1}x_{i,k}$ from $y_{i,k}$, the sign of the slopes in the transformed datasets determines whether a pair of points contributes to 
$P(\beta)$ or $Q(\beta)$ from Definition~\ref{defCtilde}.
In particular, all figures in the proof section
are plotted with transformed coordinates, i.e. $(x_{i,k}',y_{i,k}')  =  (\beta^{-1}x_{i,k},y_{i,k} - \beta^{-1}x_{i,k})$.

\subsection{Proof of Theorem \ref{T1}}

\proof{

By defining 
\begin{align*}
	Z_{ij}^{kl} := \begin{cases}
						 sgn(S_{ij}^{kl}),& \quad k \neq l \\
						 0,& \quad k = l
					\end{cases}	 
\end{align*}
as the sign of the slope in the transformed dataset, it follows that
\begin{align*}
	\tilde{C} = \sum_{\substack{k,l = 1\\k < l}}^m \sum_{i = 1}^{p_k}\sum_{j = 1}^{p_l}  Z_{ij}^{kl}.
\end{align*}
We will hence compute
\begin{align*}
	\mathbb V\bigg( \sum_{\substack{k,l = 1\\k < l}}^m \sum_{i = 1}^{p_k}\sum_{j = 1}^{p_l}  Z_{ij}^{kl} \bigg) =& ~
	\mathbb V\bigg( \frac{1}{2} \sum_{\substack{k,l = 1\\k \neq l}}^m \sum_{i = 1}^{p_k}\sum_{j = 1}^{p_l}  Z_{ij}^{kl} \bigg) \\
	=& ~ \frac{1}{4}  \bigg( \mathbb E\big( ( \sum_{\substack{k,l = 1\\k \neq l}}^m \sum_{i = 1}^{p_k}\sum_{j = 1}^{p_l}  Z_{ij}^{kl})^2  \big) - \big( \mathbb E ( \sum_{\substack{k,l = 1\\k \neq l}}^m \sum_{i = 1}^{p_k}\sum_{j = 1}^{p_l}  Z_{ij}^{kl})  \big)^2 \bigg).
\end{align*}
Since $\mathbb E(Z_{ij}^{kl}) = 0$ for all $\big( (i,k),(j,l) \big)$, the second term on the right hand side vanishes and we have 
\begin{align}
	\mathbb V\bigg( \sum_{\substack{k,l = 1\\k < l}}^m \sum_{i = 1}^{p_k}\sum_{j = 1}^{p_l}  Z_{ij}^{kl} \bigg) = \frac{1}{4} \mathbb  E \bigg( ( \sum_{\substack{k,l = 1\\k \neq l}}^m \sum_{i = 1}^{p_k}\sum_{j = 1}^{p_l}  Z_{ij}^{kl})^2 \bigg).
	\label{FormelAPBVar}
\end{align}
We split the right hand side into the sums over the squares $(Z_{ij}^{kl})^2$ and the 
sum over the remaining mixed terms $Z_{ij}^{kl} \cdot Z_{rs}^{tu}$:
\begin{align*}
	\mathbb E \bigg( ( \sum_{\substack{k,l = 1\\k \neq l}}^m \sum_{i = 1}^{p_k}\sum_{j = 1}^{p_l}  Z_{ij}^{kl})^2 \bigg) 
	= ~  \mathbb E \bigg( \sum_{\substack{k,l = 1\\k \neq l}}^m \sum_{i = 1}^{p_k}\sum_{j = 1}^{p_l}  (Z_{ij}^{kl})^2 \bigg) 
	+ \mathbb E \bigg( \sum_{\substack{k,l = 1\\k \neq l}}^m \sum_{i = 1}^{p_k}\sum_{j = 1}^{p_l}  Z_{ij}^{kl} \cdot 
	  \sum_{\substack{(r,s,t,u)\\\neq\\(i,j,k,l)}} Z_{rs}^{tu} \bigg)
\end{align*}

(a)
The first term on the right hand side is not hard to compute since 
$(Z_{ij}^{kl})^2 = 1$ for all $ \big( (i,k),(j,l) \big)$ with $k \neq l$ holds. We get
\begin{align}
  \sum_{\substack{k,l = 1\\k \neq l}}^m \sum_{i = 1}^{p_k}\sum_{j = 1}^{p_l}  (Z_{ij}^{kl})^2 
  = 2 (n^2 - \sum\limits_{k=1}^m p_k^2 )
  = 2 \big(n(n-1) - \sum_{k=1}^{m} p_k ( p_{k}-1 ) \big). \label{partA}
\end{align}

(b)
For the computation of the mixed terms, we only need to consider those $Z_{ij}^{kl} \cdot Z_{rs}^{tu}$ with a common pair of indices,
i.e.\ one of the four cases $(i,k) = (r,t), (i,k) = (s,u), (j,l)=(r,t)$ and $(j,l)=(s,u)$.
Ignoring group assignments there are $4n(n-1)(n-2)$ such combinations.
All other combinations vanish in expectation, due to independence. 
For the expectation of the sum of mixed terms, we look at the covariances of the $Z_{ij}^{kl}$ with each other.
Since all cases lead to the same result, we can assume without loss of generality  $(j,l) = (r,t)$ and multiply
the result by four.
We will combine triplets of mixed terms to simplify the calculations.
By distinction of cases as sketched in Figure \ref{KovFall2}, we obtain for the expectation
\begin{equation}
\begin{aligned}
	&\mathbb E(Z_{ij}^{kl} \cdot Z_{is}^{ku}) + \mathbb E(Z_{ij}^{kl} \cdot Z_{js}^{lu}) + \mathbb E(Z_{is}^{ku} \cdot Z_{js}^{lu})  \\
	&= \frac{1}{3}  \big( 1+1+1\big) + \frac{2}{3} \big( 1+(-1)+(-1)\big) \\
	&= \frac{1}{3} 
\end{aligned}
\label{Cov}
\end{equation}
considering that the cases occur with probability $\frac{1}{3}$ and $\frac{2}{3}$, respectively. \\

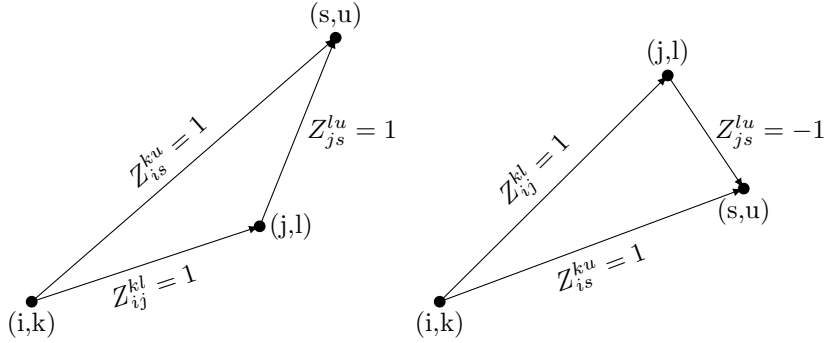
\begin{figure}[h]
	\begin{center}
		\begin{tikzpicture}
		\filldraw (0,0) circle (2pt) node[align=center, below]{(i,k)};
		\filldraw (3,1) circle (2pt) node[align=center, right]{(j,l)};
		\filldraw (4,3.5) circle (2pt) node[align=center, above]{(s,u)};
		\draw[->, >=latex] (0,0) -- (3,1) node[pos=0.5,below,sloped] {$Z_{ij}^{kl} = 1$}; 
		\draw[->, >=latex] (3,1) -- (4,3.5) node[pos=0.5,right] {$Z_{js}^{lu} = 1$};
		\draw[->, >=latex] (0,0) -- (4,3.5) node[pos=0.5,above,sloped] {$Z_{is}^{ku} = 1$};
		\end{tikzpicture}%
		\begin{tikzpicture}
		\filldraw (0,0) circle (2pt) node[align=center, below]{(i,k)};
		\filldraw (3,3) circle (2pt) node[align=center, above]{(j,l)};
		\filldraw (4,1.5) circle (2pt) node[align=center, below]{(s,u)};
		\draw[->, >=latex] (0,0) -- (3,3) node[pos=0.5,above,sloped] {$Z_{ij}^{kl} = 1$}; 
		\draw[->, >=latex] (0,0) -- (4,1.5) node[pos=0.5,below,sloped] {$Z_{is}^{ku} = 1$};
		\draw[->, >=latex] (3,3) -- (4,1.5) node[pos=0.5,right] {$Z_{js}^{lu} = -1$};
		\end{tikzpicture}
	\end{center}
	\caption{Left: Sketch of case 1. With probability $\frac13$, $y_{j,l}$ is inbetween $y_{i,k}$ and $y_{s,u}$. Slopes have equal signs. 
	Right: Sketch of case 2. With probability $\frac23$, $y_{j,l}$ is above or below both, $y_{i,k}$ and $y_{s,u}$. Slopes have different signs.} 
	\label{KovFall2}
\end{figure}

Furthermore, we need to consider that we counted each pair of slopes three times due to the arrangement in triples.
Consequently, we have to divide by $3$ and get the interim result

\begin{equation}
 4 n(n-1)(n-2) \frac19 \label{partB}
\end{equation}

So far we did not consider the groups.
Since the method from Definition~\ref{def:extendedPB} does not consider any pairs of data points within the same group, 
we have to subtract those triples of indices $\big( (i,k), (j,l), (s,u) \big)$ with $k=l, k=u$, and $l=u$, as well as $k=l=u$. \\

(c)
There are 
\begin{equation}
\sum\limits_{k=1}^m p_k (p_k -1) (n-p_k) = \sum\limits_{i=1}^m p_k (p_k -1) \sum\limits_{u\neq k}^m p_u \label{partC_count} 
\end{equation}
possible combinations of three elements from two different groups.
Let's say $k=l$. In that case, the only remaining term in~\eqref{Cov} is $Cov(Z_{is}^{ku}, Z_{js}^{lu})$ and
we have to subtract the other two terms from our computation of \eqref{FormelAPBVar} because we wrongfully included them in step (b).
In case of strictly separated groups, however, those combinations vanish as can be seen below in equation \eqref{partC}. 

In case of two equal indices 
there are three different possibilities, how the points are located relativ to each other. 
Looking at the single point from group $u$, it can be located above, in between or below the two points from group $k$ (relating to the y-axis).
Those arrangements occur with probability $\frac{1}{3}$ each and are illustrated in Figure~\ref{3el2groups} if $x_{s,u} > \max(x_{i,k},x_{j,k})$ 
and in Figure~\ref{3el2groupsinbetween} for $\min(x_{i,k},x_{j,k}) < x_{s,u} < \max(x_{i,k},x_{j,k})$.
Summed up, those arrangements yield an expectation of 
\begin{equation}
	\begin{aligned}
		&\qquad \mathbb E(Z_{ij}^{kk} \cdot Z_{is}^{ku}) + \mathbb E(Z_{ij}^{kk} \cdot Z_{js}^{ku}) \\
		&= (1-(q_{ku})) \frac13 \big( 0 + 2 -2 \big) + q_{ku} \frac13 \big( 2 + 0 +0 \big)
                = \frac 23 q_{ku}
	\end{aligned}
	\label{partC}
\end{equation}

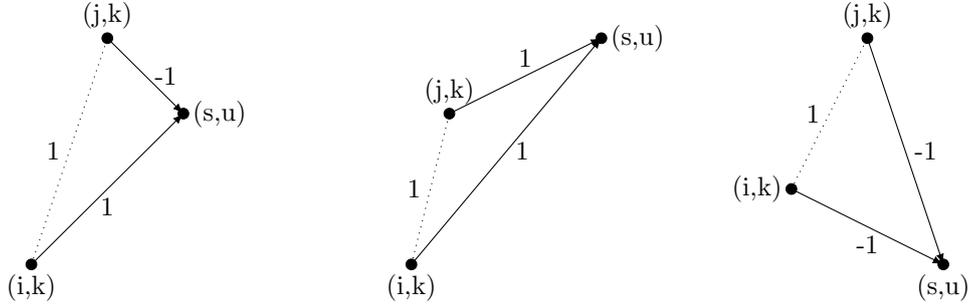
\begin{figure}[h]
	\begin{center}
			\begin{tikzpicture}
				\filldraw (0,0) circle (2pt) node[align=center, below]{(i,k)};
				\filldraw (1,3) circle (2pt) node[align=center, above]{(j,k)};
				\filldraw (2,2) circle (2pt) node[align=center, right]{(s,u)};
				\draw[dotted, >=latex] (0,0) -- (1,3) node[pos=0.5,left] {1}; 
				\draw[->, >=latex] (0,0) -- (2,2) node[pos=0.5,below] {1};
				\draw[->, >=latex] (1,3) -- (2,2) node[pos=0.5,right] {-1};
				\filldraw (5,0) circle (2pt) node[align=center, below]{(i,k)};
				\filldraw (5.5,2) circle (2pt) node[align=center, above]{(j,k)};
				\filldraw (7.5,3) circle (2pt) node[align=center, right]{(s,u)};
				\draw[dotted, >=latex] (5,0) -- (5.5,2) node[pos=0.5,left] {1}; 
				\draw[->, >=latex] (5,0) -- (7.5,3) node[pos=0.5,right] {1};
				\draw[->, >=latex] (5.5,2) -- (7.5,3) node[pos=0.5,above] {1};	
				\filldraw (10,1) circle (2pt) node[align=center, left]{(i,k)};
				\filldraw (11,3) circle (2pt) node[align=center, above]{(j,k)};
				\filldraw (12,0) circle (2pt) node[align=center, below]{(s,u)};
				\draw[dotted, >=latex] (10,1) -- (11,3) node[pos=0.5,left] {1}; 
				\draw[->, >=latex] (10,1) -- (12,0) node[pos=0.5,below] {-1};
				\draw[->, >=latex] (11,3) -- (12,0) node[pos=0.5,right] {-1};							
			\end{tikzpicture}
	\end{center}
	\caption{All possible arrangements of three elements from two different groups. 
	$(i,k)$ and $(j,k)$ come from the same group, $(s,u)$ comes from another group and $x_{s,u} > \max(x_{i,k},x_{j,k})$.}
	\label{3el2groups}
\end{figure}

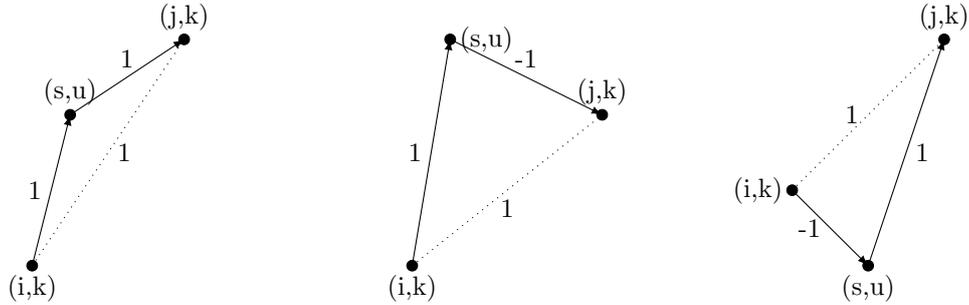
\begin{figure}[h]
	\begin{center}
			\begin{tikzpicture}
				\filldraw (0,0) circle (2pt) node[align=center, below]{(i,k)};
				\filldraw (2,3) circle (2pt) node[align=center, above]{(j,k)};
				\filldraw (0.5,2) circle (2pt) node[align=center, above]{(s,u)};
				\draw[dotted, >=latex] (0,0) -- (2,3) node[pos=0.5,right] {1}; 
				\draw[->, >=latex] (0,0) -- (0.5,2) node[pos=0.5,left] {1};
				\draw[<-, >=latex] (2,3) -- (0.5,2) node[pos=0.5,above] {1};
				\filldraw (5,0) circle (2pt) node[align=center, below]{(i,k)};
				\filldraw (7.5,2) circle (2pt) node[align=center, above]{(j,k)};
				\filldraw (5.5,3) circle (2pt) node[align=center, right]{(s,u)};
				\draw[dotted, >=latex] (5,0) -- (7.5,2) node[pos=0.5,below] {1}; 
				\draw[->, >=latex] (5,0) -- (5.5,3) node[pos=0.5,left] {1};
				\draw[<-, >=latex] (7.5,2) -- (5.5,3) node[pos=0.5,above] {-1};	
				\filldraw (10,1) circle (2pt) node[align=center, left]{(i,k)};
				\filldraw (12,3) circle (2pt) node[align=center, above]{(j,k)};
				\filldraw (11,0) circle (2pt) node[align=center, below]{(s,u)};
				\draw[dotted, >=latex] (10,1) -- (12,3) node[pos=0.5,left] {1}; 
				\draw[->, >=latex] (10,1) -- (11,0) node[pos=0.5,left] {-1};
				\draw[<-, >=latex] (12,3) -- (11,0) node[pos=0.5,right] {1};							
			\end{tikzpicture}
	\end{center}
	\caption{Some as Figure~\ref{3el2groups} if $\min(x_{i,k},x_{j,k}) < x_{s,u} < \max(x_{i,k},x_{j,k})$ with probability $q_{ku}$.}
	\label{3el2groupsinbetween}
\end{figure}

(d)
Further, in \eqref{partB} we wrongfully included the case $k=l=u$. This means, we have to substract
\begin{equation}
  \frac 49 \sum\limits_{k=1}^m p_k (p_k-1)(p_k-2) \label{partD}
\end{equation}
and can combine \eqref{partA} and \eqref{partB}-
\eqref{partD} to get for the term in \eqref{FormelAPBVar}
\begin{align}
 \frac12 \big( n(n-1) - \sum\limits_{k=1}^m p_k(p_k-1) \big)    +  \frac19 \big( n(n-1)(n-2)  -  \sum\limits_{k=1}^m p_k(p_k-1)(p_k-2 + 2 \sum\limits_{u\neq k}^m p_u q_{ku} ) \big)  \nonumber
\end{align}

the wanted formula for the variance. Naturally, this formula holds likewise for the less complex settings
of equally sized non-overlapping groups with $p_k = p$ for all $k$, i.e.
\begin{align}
 \mathbb V[\tilde{C}] &= \frac12 n(n-p) +  \frac{n}{9} \big((n-1)(n-2)  -  (p-1)(p-2) \big)\\
		      &= \frac{n}{18} \big( 3(n-p) + 2(n^2-p^2) \big)
\end{align}
and for non-grouped data with $m = n$ and $p_k = 1$ we regain the classic result
\begin{equation}
 \mathbb V[C] = \frac{n(n-1)(2n+5)}{18}.
\end{equation}

}

\subsection{Proof of Theorem \ref{T2}}

\proof{

  The proofs for the PBR and TSR for non-grouped measurements \cite{Passing1983,Sen1968} relied on results for rank correlation methods.
  Our proof will follow the same route as shown in chapter 5 of Kendalls book~\cite{Kendall1990rank}.
  To show the asymptotic normality we will show that the moments 
  of the distribution of $\tilde{C}$ tend to those of the normal distribution~\cite[Thm~30.2]{Billingsley1995}.

  Under the null hypothesis we have that $\mathbb P(Z^{i,k}_{j,l} = 1) = \mathbb P(Z^{i,k}_{j,l} = -1) $.
  Hence moments of $\tilde{C}$ of odd order vanish.
  For moments of even order, we need to compute
  
  \begin{equation}
  \mathbb E[ \big( \sum\limits \sum\limits  Z_{ij}^{kl} \big)^{2r} ]. \label{evenmoments}
  \end{equation}
  
  Consider the expansion of the sum above, which consists of 
  summands with $2r$ factors each.
  Since $\mathbb E [ Z_{ij}^{kl} ] = 0$ and the  $Z_{ij}^{kl}$ are independent if $k \neq l$
  each summand with an independent factor will vanish in the expectation.
  Let us consider the number of summands where the factors are pairwise linked by exactly one suffix.
  Each of these summands will look like this:
  \begin{equation}
    \mathbb E[ Z_{ij}^{kl} Z_{is}^{ku} Z_{rv}^{tw} Z_{rx}^{ty} \dots  ], \label{summand}
  \end{equation}
  which due to independence could be split up into the form
  $$ (\mathbb E[Z_{ij}^{kl} Z_{is}^{ku}])^r . $$
  
  However, by doing so we would no longer connect the cases as in equation \eqref{Cov}. So we ignore the above and combine cases,
  as we have done in the proof of Theorem~\ref{T1}:
  
  For each $Z_{ij}^{kl} Z_{is}^{ku}$ in each summand as shown in equation~\eqref{summand} 
  there are two other combinations that are exactly the same except that there is an
  $Z_{ij}^{kl} Z_{jl}^{us}$ or $Z_{is}^{ku} Z_{jl}^{us}$ instead of $Z_{ij}^{kl} Z_{is}^{ku}$, respectively.
  Therefore we define for each triplet of points $t = ((i,k),(j,l),(s,u))$ the pairings
  
  \begin{align*}
   p_1((i,k),(j,l),(s,u)) &= Z_{ij}^{kl} Z_{is}^{ku}\\
   p_2((i,k),(j,l),(s,u)) &= Z_{ij}^{kl} Z_{js}^{lu}\\
   p_3((i,k),(j,l),(s,u)) &= Z_{is}^{ku} Z_{js}^{lu}\\
  \end{align*}

  For each summand with pairwise tied indices let $t_1,\dots,t_r$ be the triplet of indices used in the $r$-th pairing.
  Let us then consider the set of $3^r$ sums which is represented by
  $ \prod\limits_{i = 1}^r \sum\limits_{j=1}^3 p_j(t_i).$  
  Now we see that  
  \begin{align*}
   \mathbb E  \Big[ \prod\limits_{i = 1}^r \sum\limits_{j=1}^3 p_j(t_i) \Big]
   = \prod\limits_{i = 1}^r \sum\limits_{j=1}^3 \mathbb E  \Big[ p_j(t_i) \Big]
   = \prod\limits_{i = 1}^r \frac 13 = \frac{1}{3^r},
  \end{align*}
  if the points of $t_1,\dots,t_r$ are a disjoint and hence independent set of $3r$ points.
  
 \pagebreak[1]
  
  The remaining questions are:
  \begin{enumerate}
   \item[(a)] How many summands with pairwise tied indices are there in the expanded sum ignoring group associations?
   \item[(b)] Which of these need to be omitted due to group associations for non-overlapping groups?
   \item[(c)] What changes for overlapping groups?
  \end{enumerate}
  
  Answer to (a):
  
  There are $\binom{2r}{r}$ ways to choose the first factors of $r$ pairings and $r!$ ways to assign the remaining $r$ factors.
  But now we have counted some combinations twice so we have to divide by $2^r$ to get the final number of ways to choose $r$ pairings.
  Given such a combination there are now $3r$ different indices to choose, so $\frac{n!}{(n-3r)!} \sim n^{3r}$ possibilities.
  
  Thus we end up with 
  $$ \binom{2r}{r} \frac{r!}{2^r} n^{3r} = \frac{(2r)!}{r!  2^r} n^{3r}$$
  combinations with $r$ pairs with tied indices.
  Finally since we counted each summand $3^r$ times we would get for~\eqref{evenmoments}
  \begin{equation}
   (2r-1)!! (\frac 19 n^3)^r  \label{classiccase}
  \end{equation}
   if we ignore group associations.
  
  Answer to (b):
  
  For $t= ((i,k),(j,l),(s,u))$ there are 3 different scenarios.
  
  \begin{description}
   \item[Scenario 1: all groups are different ($k\neq l\neq u$)] \quad\\ No difference to the above calculation
   \item[Scenario 2: all groups are equal ($k=l=u$)] \quad\\The corresponding summand is not included in~\eqref{evenmoments}
   which we mimic by setting\\ $\sum \limits_{j=1}^3 \mathbb E  \Big[p_j(t)\Big]=0$
   \item[Scenario 3: two groups are equal ($k=l\neq u$ or $ k \neq l= u $ or $ k=u\neq l $)] \quad\\
   Althoug some summands are not included in this case, we do not have to change anything
   since the effect vanishes for separated groups as shown in equation \eqref{partC}.
  \end{description}

  For separated groups only scenario 2 alters the above calculation. 
  So we have to subtract the number of combinations with at
  least one factor where all indices are from the same group.
  The probablity to choose 3 indices from group $k$ is given by 
  $n^{-3} p_k(n) ( p_k(n) -1) (p_k(n) -2) $ and thus the 
  fraction of all sums that have no $t_i$ from scenario 2 tends to
  \begin{equation}
   a^r := \bigg( 1 - \lim\limits_{n\to\infty} \sum\limits_{k=1}^{m} \frac{\big(p_k(n)\big)^3} {n^3} \bigg)^r, \label{fraction}
  \end{equation}
  which simplifies for $m$ groups with equal sizes to
  \begin{equation}
  \big(1 - \frac{1}{m^2} \big)^r. 
  \end{equation}

  So we have to multiply equation \eqref{classiccase} and \eqref{fraction} and end up with
  $$ \mathbb E[ \big( \sum\limits \sum\limits  Z_{ij}^{kl} \big)^{2r} ] = (2r-1)!!
  \bigg(\frac 19 n^3 \big(1 - \lim\limits_{n\to\infty} \sum\limits_{k=1}^{m} \frac{\big(p_k(n)\big)^3} {n^3} \big) \bigg)^r.$$

  Answer to (c):
  
  If groups can overlap, let $r_1$ be the number of
  triplets where groups do not overlap (Scenario 3.1) and $r_2$ the number of triplets where groups overlap (Scenario 3.2).
  
  In this case we get
  \begin{align*}
   \mathbb E  \Big[ \prod\limits_{i = 1}^r \sum\limits_{j=1}^3 p_j(t_i) \Big]
   = (-1)^{r_2} \frac{1}{3^r}.
  \end{align*}  
  So for each summand with an odd number of triplets with $\min(x_{i,k},x_{j,k}) < x_{s,u} < \max(x_{i,k},x_{j,k})$ we have to substract
  $\frac{2}{3^r}$ just as in the proof of Theorem~\ref{T1}.
  We set $$b = \lim\limits_{n\to \infty} \sum\limits_{k=1}^m \sum\limits_{u\neq k}^m n^{-3} \big(p_k(n)\big)^2 p_u(n) q_{ku}  $$
  The probability to get a summand with odd $r_2$ given there is no triplet from scenario 2 is now given by 
  $ \frac 12 (1-(1-2 \frac{b}{a})^r).$
  This leads us to 
  \begin{align*}
  &\frac{n^{3r}}{9^r} \left(a^r - 2 a^r ( \frac 12 (1-(1-2 \frac{b}{a})^r) )  \right)  \\
  = &\frac{n^{3r}}{9^r} \left(a^r - a^r   (1-  (1 +  \sum\limits_{i=1}^r \binom{r}{i} \big(\frac{-2b}{a}\big)^i  )       )   \right)  \\
  = &\frac{n^{3r}}{9^r} \left(a^r +   \sum\limits_{i=1}^r \binom{r}{i} \big(-2b\big)^i a^{r-i}           \right)
  = \frac{n^{3r}}{9^r} \left(a^r - 2b  \right)^r.
  \end{align*}
  As all other terms are of lower order in $n$ this completes the proof.
  To see this consider the number of summands where the factors are all linked but not pairwise, 
  i.e.\ there is at least one summand that is linked to more than one other summand.
  In this case, there are less than $3r-1$ indices to choose, such that the order is not greater than $n^{3r-1}$.
  
}

\subsection{Justification of confidence intervals in Corollary~\ref{T3}}

For the justification of the formula $M_1 = \floor{ \frac{N - \tilde C_\gamma}{2} }$ 
we proceed just as Passing and Bablok~\cite{Passing1983}. Consider
\begin{align*}
	S_{(M_1 + K)} < \beta < S_{(M_2 + K)}
\end{align*}
which holds if and only if $M_1 + K \leq ~ Q(\beta)$ and $M_1 -  K \leq ~ P(\beta)$
hold and thus if and only if
\begin{align*}
	2M_1 - N \leq P(\beta) - Q(\beta) +2K \leq N - 2M_1.
\end{align*}
We know that the distribution of $\tilde C:= P(\beta) - Q(\beta)$ does not depend on the distribution of $(X,Y)$, 
whereas the distribution of $K$ does.
Therefore, we can not provide a general formula for $M_1$ satisfying
\begin{align*}
	P \big( S_{(M_1 + K)} < \beta < S_{(M_2 + K)}  \big) = &~ P \big( 2M_1 - N \leq \tilde C +2K \leq N - 2M_1  \big) \\
	=& ~ 1 - \alpha
\end{align*}
that does not dependent on the distribution of $(X,Y)$. 
However, 
as shown in Theorem~\ref{T2}, $\tilde C$ is asymptotically normally distributed, 
such that $P\big( -\tilde C_\gamma \leq \tilde C \leq \tilde C_\gamma \big) \sim 1 - \gamma$. Thus, we argue just as Passing and Bablok~\cite{Passing1983} that
$M_1$ can be found by setting 
	$N - 2M_1 = \tilde C_\gamma$.

\section*{acknowlegement}

We thank Peter Pfaffelhuber for various discussions pointing us in the right direction as well as
Lukas Steinberger and an anonymous reviewer for helpful comments and careful reading of the manuscript.

\newpage

\bibliographystyle{abbrv}
\bibliography{passingbablok.bib}

\end{document}